\documentclass{amsproc}
\usepackage{amssymb}
\usepackage{graphicx}
\usepackage{amscd}
\usepackage{amsmath}

\newtheorem{definition}{Definition}

\newtheorem{proposition}{Proposition}
\newtheorem{remark}{Remark}

\newtheorem{theorem}{Theorem}

\newtheorem{agreement}{Agreement}
\numberwithin{equation}{section}

\begin{document}

\title{On an application of the Boundary control method to classical
moment problems.}

\author{A S  Mikhaylov$^1,^2$, V S  Mikhaylov$^1$} 
\address{$^1$St. Petersburg   Department   of   V.A. Steklov    Institute   of   Mathematics
of   the   Russian   Academy   of   Sciences, 7, Fontanka, 191023
St. Petersburg, Russia}
\address{$^2$Saint Petersburg State University, St.Petersburg State
University, 7/9 Universitetskaya nab., St. Petersburg, 199034
Russia.}



\maketitle

\begin{abstract} We establish relationships between the classical moments problems
which are problems of a construction of a measure supported on a
real line, on a half-line or on an interval from prescribed set of
moments with the Boundary control approach to a dynamic inverse
problem for a dynamical system with discrete time associated with
Jacobi matrices. We show that the solution of corresponding
truncated moment problems is equivalent to solving some
generalized spectral problems.

\end{abstract}

\section{Introduction}

In \cite{MM4} the authors put forward an approach to Hamburger,
Stieltjes and Hausdorff moment problems based on their
relationships with inverse problems for dynamical systems with
discrete time associated with Jacobi matrices. In the present
paper we utilize some ideas from \cite{MikhaylovAS_MM2} about de
Branges spaces associated with such dynamical systems and extend
and elaborate results obtained in \cite{MM4}. We begin with
introducing moment problems and spaces of polynomials associated
with it, dynamical systems with discrete time associated with
Jacobi matrices and de Branges spaces of analytic functions.

\subsection{Classical moment problems.}

For a given a sequence of numbers $s_0,s_1,s_2,\ldots$ called
moments, a solution of a Hamburger moment problem
\cite{MikhaylovAS_Ahiez,MikhaylovAS_S} is a Borel measure
$d\rho(\lambda)$ on $\mathbb{R}$ such that
\begin{equation}
\label{MikhaylovAS_Moment_eq} s_k=\int_{-\infty}^\infty
\lambda^k\,d\rho(\lambda),\quad k=0,1,2,\ldots
\end{equation}
The measure is a solution to Stieltjes or Hausdorff moment
problems provided $\operatorname{supp}d\rho\subset (0,+\infty)$ or
$\operatorname{supp}d\rho\subset (0,1)$ respectively; in these
cases the moments are called Hamburger, Stieltjes or Hausdorff.

Following \cite{MikhaylovAS_S,Schm} we denote by $C[X]$ the set of
complex polynomials and by $C_N[X]$ the set of polynomials of
order less than or equal to $N$. The moments
$\left\{s_k\right\}_{k=0}^\infty$ determine on $C[X]$ the bilinear
form by the rule: for $F,G\in C[X]$,
$F(\lambda)=\sum_{n=0}^{N-1}\alpha_n \lambda^n,$
$G(\lambda)=\sum_{n=0}^{N-1}\beta_n \lambda^n,$ one defines
\begin{equation}
\label{Moment_scal} \langle
F,G\rangle=\sum_{n,m=0}^{N-1}s_{n+m}\alpha_n\overline{\beta_m}.
\end{equation}
Thus this quadratic form is determined by the following
(semi-infinite) Hankel matrix:
\begin{equation}
\label{S_matr}
S=\begin{pmatrix} s_0 & s_1 & s_2 & s_3 & \ldots\\
s_1 & s_2 & s_3 & \ldots & \ldots \\
s_2 & s_3 & \ldots & \ldots & \ldots \\
s_3 & \ldots & \ldots & \ldots & \ldots \\
\ldots & \ldots & \ldots & \ldots & \ldots
\end{pmatrix}
\end{equation}

\subsection{Initial boundary value problems associated with
Jacobi matices.}

An initial boundary value problem (IBVP) for an auxiliary
dynamical system with discrete time for a Jacobi matrix is set up
in the following way: for a given sequence of positive numbers
$\{a_0, a_1,\ldots\}$ (in what follows we assume $a_0=1$) and real
numbers $\{b_1, b_2,\ldots \}$, we denote by $A$ the Jacobi
operator, defined on $l_2$, which has a matrix form:
\begin{equation}
\label{MikhaylovAS_Jac_matr}
A=\begin{pmatrix} b_1 & a_1 & 0 & 0 & 0 &\ldots \\
a_1 & b_2 & a_2 & 0 & 0 &\ldots \\
0 & a_2 & b_3 & a_3 & 0 & \ldots \\
\ldots &\ldots  &\ldots &\ldots & \ldots &\ldots
\end{pmatrix}
\end{equation}
All our considerations will be local, so when coefficients
$a_j,b_j$ are such that Jacobi matrix is in limit circle case, by
$A$ we can assume any self-adjoint extension. For $N\in
\mathbb{N}$, by $A^N$ we denote the $N\times N$ Jacobi matrix
which is a block of (\ref{MikhaylovAS_Jac_matr}) consisting of the
intersection of first $N$ columns with first $N$ rows of $A$. We
consider the dynamical system with discrete time associated with
$A^N$:
\begin{equation}
\label{MikhaylovAS_Jacobi_dyn_int}
\begin{cases}
v_{n,t+1}+v_{n,t-1}-a_nv_{n+1,t}-a_{n-1}v_{n-1,t}-b_nv_{n,t}=0,\,\, t\in \mathbb{N}\cup\{0\},\,\, n\in 1,\ldots, N,\\
v_{n,\,-1}=v_{n,\,0}=0,\quad n=1,2,\ldots,N+1, \\
v_{0,\,t}=f_t,\quad v_{N+1,\,t}=0,\quad t\in \mathbb{N}_0,
\end{cases}
\end{equation}
where $f=(f_0,f_1,\ldots)$ is a \emph{boundary control}. The
solution to (\ref{MikhaylovAS_Jacobi_dyn_int}) is denoted by
$v^f$. Note that (\ref{MikhaylovAS_Jacobi_dyn_int}) is a discrete
analog of dynamical system with boundary control for a wave
equation on an interval \cite{MikhaylovAS_AM,MikhaylovAS_BM_1}.

The operator corresponding to a finite Jacobi matrix we also
denote by $A^N: \mathbb{R}^N\mapsto \mathrm{R}^N$, it is given by
\begin{equation}
\label{A_eq1}
\begin{cases}
(A\psi)_n&=a_{n}\psi_{n+1}+a_{n-1}\psi_{n-1}+b_n\psi_n,\quad
2\leqslant n\leqslant N-1,\\
(A\psi)_1&=b_1\psi_1+a_1\psi_2,\quad n=1,
\end{cases}
\end{equation}
and the Dirichlet condition at the "right end":
\begin{equation}
\label{A_eq2} \psi_{N+1}=0.
\end{equation}
We also consider the dynamical system corresponding to a
semi-infinite Jacobi matrix:
\begin{equation}
\label{Jacobi_dyn}
\begin{cases}
u_{n,t+1}+u_{n,t-1}-a_nu_{n+1,t}-a_{n-1}u_{n-1,t}-b_nu_{n,t}=0,\,\, t\in \mathbb{N}\cup\{0\},\,\, n\in 1,\ldots, N,\\
u_{n,\,-1}=u_{n,\,0}=0,\quad n=1,2,\ldots,N+1, \\
u_{0,\,t}=f_t,\quad t\in \mathbb{N}_0,
\end{cases}
\end{equation}
its solution is denoted by $u^f.$

\subsection{De Branges spaces.}

Here we provide the information on de Branges spaces in accordance
with \cite{RR}. The entire function $E:\mathbb{C}\mapsto
\mathbb{C}$ is called a \emph{Hermite-Biehler function} if
$|E(z)|>|E(\overline z)|$ for $z\in \mathbb{C}_+$. We use the
notation $F^\#(z)=\overline{F(\overline{z})}$. The \emph{Hardy
space} $H_2$ is defined by: $f\in H_2$ if $f$ is holomorphic in
$\mathbb{C}^+$ and
$\sup_{y>0}\int_{-\infty}^\infty|f(x+iy)|^2\,dx<\infty$. Then the
\emph{de Branges space} $B(E)$ consists of entire functions such
that:
\begin{equation*}
B(E):=\left\{F:\mathbb{C}\mapsto \mathbb{C},\,F \text{ entire},
\int_{\mathbb{R}}\left|\frac{F(\lambda)}{E(\lambda)}\right|^2\,d\lambda<\infty,\,\frac{F}{E},\frac{F^\#}{E}\in
H_2\right\}.
\end{equation*}
The space $B(E)$ with the scalar product
\begin{equation*}
[F,G]_{B(E)}=\frac{1}{\pi}\int_{\mathbb{R}}{F(\lambda)}
\overline{G(\lambda)}\frac{d\lambda}{|E(\lambda)|^2}
\end{equation*}
is a Hilbert space. For any $z\in \mathbb{C}$ the
\emph{reproducing kernel} is introduced by the relation
\begin{equation}
\label{repr_ker} J_z(\xi):=\frac{\overline{E(z)}E(\xi)-E(\overline
z)\overline{E(\overline \xi)}}{2i(\overline z-\xi)}.
\end{equation}
Then
\begin{equation*}
F(z)=[J_z,F]_{B(E)}=\frac{1}{\pi}\int_{\mathbb{R}}{J_z(\lambda)}
\overline{G(\lambda)}\frac{d\lambda}{|E(\lambda)|^2}.
\end{equation*}
We observe that a Hermite-Biehler function $E(\lambda)$ defines
$J_z$ by (\ref{repr_ker}).

In the second section we provide the results on solutions to
(\ref{MikhaylovAS_Jacobi_dyn_int}) and (\ref{Jacobi_dyn}) and
introduce the operators of the Boundary control method according
to \cite{MikhaylovAS_MM3}; using the ideas of
\cite{MikhaylovAS_MM2} we introduce de Branges spaces
corresponding to dynamical systems (\ref{Jacobi_dyn}),
(\ref{MikhaylovAS_Jacobi_dyn_int}) and give representation of
reproducing kernel in the space of polynomials and Christoffel
symbols \cite{Schm} in dynamic terms. In the third section we
outline the solution to a truncated moment problem following
\cite{MM4}, more specifically, we reduce it to the generalized
spectral problem for special matrices constructed from moments. In
the last section we apply obtained results to the problem of
uniqueness of solutions to moment problems.

\section{IBVP for a dynamical system associated with Jacobi matrix and de Branges spaces.}

We fix some positive integer $T$ and denote by $\mathcal{F}^T$ the
\emph{outer space} of systems (\ref{MikhaylovAS_Jacobi_dyn_int}),
(\ref{Jacobi_dyn}), the space of controls:
$\mathcal{F}^T:=\mathbb{R}^T$, $f\in \mathcal{F}^T$,
$f=(f_0,\ldots,f_{T-1})$, we use the notation
$\mathcal{F}^\infty=\mathbb{R}^\infty$ when control acts for all
$t\geqslant 0$.
\begin{definition}
For $f,g\in \mathcal{F}^\infty$ we define the convolution
$c=f*g\in \mathcal{F}^\infty$ by the formula
\begin{equation*}
c_t=\sum_{s=0}^{t}f_sg_{t-s},\quad t\in \mathbb{N}\cup \{0\}.
\end{equation*}
\end{definition}

The input $\longmapsto$ output correspondences in systems
(\ref{MikhaylovAS_Jacobi_dyn_int}), (\ref{Jacobi_dyn}) are
realized by a \emph{response operators}:
$R^T_N,\,R^T:\mathcal{F}^T\mapsto \mathbb{R}^T$ defined by rules
\begin{eqnarray*}
\left(R^T_Nf\right)_t=v^f_{1,\,t}=\left(r^N*f_{\cdot-1}\right)_t,
\quad t=1,\ldots,T,\\
\left(R^Tf\right)_t=u^f_{1,\,t}=\left(r*f_{\cdot-1}\right)_t,
\quad t=1,\ldots,T,
\end{eqnarray*}
where $r^N=(r_0^N,r_1^N,\ldots,r_{T-1}^N)$,
$r=(r_0,r_1,\ldots,r_{T-1})$ are \emph{response vectors},
convolution kernels of response operators. These operators play
the role of dynamic inverse data, corresponding inverse problems
were considered in \cite{MikhaylovAS_MM,MikhaylovAS_MM3}. By
choosing the special control
 $f=\delta:=(1,0,0,\ldots)$, kernels of response operators can be determined as
\begin{equation*}
\left(R^T_N\delta\right)_t=v^\delta_{1,\,t}= r_{t-1}^N,\quad
\left(R^T\delta\right)_t=u^\delta_{1,\,t}= r_{t-1}.
\end{equation*}
Let $\phi_n(\lambda)$ be a solution to the following difference
equation
\begin{equation}
\label{MikhaylovAS_Phi_def} \begin{cases} a_n\phi_{n+1}+a_{n-1}\phi_{n-1}+b_n\phi_n=\lambda\phi_n,\\
\phi_0=0,\,\,\phi_1=1.
\end{cases}
\end{equation}
Denote by $\{\lambda_k\}_{k=1}^N$ roots of the equation
$\phi_{N+1}(\lambda)=0$, it is known
\cite{MikhaylovAS_Ahiez,MikhaylovAS_S} that they are real and
distinct. We introduce vectors $\phi^n\in \mathbb{R}^N$ by the
rule $\phi^n_i:=\phi_i(\lambda_n)$, $n,i=1,\ldots,N,$ and define
numbers $\rho_k$ by
\begin{equation*}
(\phi^k,\phi^l)=\delta_{kl}\rho_k,\quad k,l=1,\ldots,N,
\end{equation*}
where $(\cdot,\cdot)$ is a scalar product in $\mathbb{R}^N$.
\begin{definition}
The set of pairs
\begin{equation*}
\{\lambda_k,\rho_k\}_{k=1}^N
\end{equation*}
is called spectral data of the operator $A^N$.
\end{definition}
Let $\mathcal{T}_k(2\lambda)$ be Chebyshev polynomials of the
second kind, i.e. they satisfy
\begin{equation}
\label{Cheb}
\begin{cases}
\mathcal{T}_{t+1}+\mathcal{T}_{t-1}-\lambda \mathcal{T}_{t}=0,\\
\mathcal{T}_{0}=0,\,\, \mathcal{T}_1=1.
\end{cases}
\end{equation}
The spectral function of $A^N$ is introduced by the rule
\begin{equation}
\label{MikhaylovAS_measure}
\rho^{N}(\lambda)=\sum_{\{k\,|\,\lambda_k<\lambda\}}\frac{1}{\rho_k},
\end{equation}
The spectral function of $A$ (non unique if $A$ is limit circle at
infinity) is denoted by $\rho(\lambda)$. In
\cite{MikhaylovAS_MM,MikhaylovAS_MM3} by the application of
Fourier expansion method following representations for the
solution $v^f$ and components of response vector were established:
\begin{proposition}
The solution to (\ref{MikhaylovAS_Jacobi_dyn_int}) and the kernel
of $R^T_N$ admit representations
\begin{eqnarray}
v^f_{n,t}=\int_{-\infty}^\infty \sum_{k=1}^t
\mathcal{T}_k(\lambda)f_{t-k}\phi_n(\lambda)\,d\rho^{N}(\lambda), \label{MikhaylovAS_Jac_sol_spectr}\\
r_{t-1}^{N}=\int_{-\infty}^\infty
\mathcal{T}_t(\lambda)\,d\rho^{N}(\lambda),\quad t\in \mathbb{N}.
\label{MikhaylovAS_Jac_resp_spectr}
\end{eqnarray}
\end{proposition}
\begin{remark}
The solution corresponding to semi-infinite Jacobi matrix $u^f$
and entries of the kernel of $R^T$ admit representations
(\ref{MikhaylovAS_Jac_sol_spectr}),
(\ref{MikhaylovAS_Jac_resp_spectr}) with $d\rho^{N}$ substituted
by $d\rho(\lambda)$.
\end{remark}

The \emph{inner space} of dynamical system
(\ref{MikhaylovAS_Jacobi_dyn_int}) is
$\mathcal{H}^N:=\mathbb{R}^N$, $h\in \mathcal{H}^N$,
$h=(h_1,\ldots, h_N)$, $v^f_{\cdot,\,T}\in \mathcal{H}^N$ for all
$T$. For the system (\ref{MikhaylovAS_Jacobi_dyn_int}) the
\emph{control operator} $W^T_{N}:\mathcal{F}^T\mapsto
\mathcal{H}^N$ is defined by the rule
\begin{equation*}
W^T_{N}f:=v^f_{n,\,T},\quad n=1,\ldots,N.
\end{equation*}
The set
\begin{equation*}
\mathcal{U}^T:=W_N^T \mathcal{F}^T=\left\{v^f_{\cdot,T}\bigl| f\in
\mathcal{F}^T\right\}
\end{equation*}
is called reachable. For the system (\ref{Jacobi_dyn}) we have
that $v^f_{\cdot,\,T}\in \mathcal{H}^T$, thus the control operator
$W^T:\mathcal{F}^T\mapsto \mathcal{H}^T$ is introduced by
\begin{equation*}
W^Tf:=u^f_{n,\,T},\quad n=1,\ldots,T.
\end{equation*}

Everywhere below we substantially use the finiteness of the speed
of wave propagation in systems (\ref{MikhaylovAS_Jacobi_dyn_int}),
(\ref{Jacobi_dyn}) which implies the following dependence of
inverse data on coefficients $\{a_n,b_n\}$: for $M\in \mathbb{N}$,
$M\leqslant N,$ the element $v^f_{1,2M-1}$ depends on
$\left\{a_1,\ldots,a_{M-1}\right\},$
$\left\{b_1,\ldots,b_{M}\right\}$. On observing this we can
formulate the following
\begin{remark}
\label{Rem1} Entries of the response vector
$(r_0^N,r_1^N,\ldots,r_{2N-2}^N)$) depend on
$\{a_0,\ldots,a_{N-1}\}$, $\{b_1,\ldots,b_{N}\}$, and does not
depend on the boundary condition at $n=N+1$, the entries starting
from $r_{2N-1}^N$ does "feel" the boundary condition at $n=N+1$.
Moreover,
\begin{equation}
\label{Rav_bc} u^f_{n,\,t}=v^f_{n,\,t},\quad n\leqslant t\leqslant
N,\quad \text{and}\quad W^N=W^N_{N}.
\end{equation}
\end{remark}

The \emph{connecting operator} for the system
(\ref{MikhaylovAS_Jacobi_dyn_int}) $C^T_{N}: \mathcal{F}^T\mapsto
\mathcal{F}^T$ is defined via the quadratic form: for arbitrary
$f,g\in \mathcal{F}^T$ we set
\begin{equation*}
\left(C^T_{N} f,g\right)_{\mathcal{F}^T}=\left(v^f_{\cdot,\,T},
v^g_{\cdot,\,T}\right)_{\mathcal{H}^N}=\left(W^T_{N}f,W^T_{N}g\right)_{\mathcal{H}^N}.
\end{equation*}
For the system (\ref{Jacobi_dyn}) the connecting operator $C^T:
\mathcal{F}^T\mapsto \mathcal{F}^T$ is introduced by the rule:
\begin{equation*}
\left(C^T f,g\right)_{\mathcal{F}^T}=\left(u^f_{\cdot,\,T},
u^g_{\cdot,\,T}\right)_{\mathcal{H}^T}=\left(W^Tf,W^Tg\right)_{\mathcal{H}^N}.
\end{equation*}

In \cite{MikhaylovAS_MM,MikhaylovAS_MM3} the following formulas
were obtained:
\begin{proposition}
Connecting operators for systems
(\ref{MikhaylovAS_Jacobi_dyn_int}), (\ref{Jacobi_dyn}) admit
spectral representations
\begin{eqnarray}
\{C^T_{N}\}_{l+1,\,m+1}=\int_{-\infty}^\infty
\mathcal{T}_{T-l}(\lambda)\mathcal{T}_{T-m}(\lambda)\,d\rho^{N}(\lambda),
\quad l,m=0,\ldots,T-1, \label{MikhaylovAS_SP_mes_d}\\
\{C^T\}_{l+1,\,m+1}=\int_{-\infty}^\infty
\mathcal{T}_{T-l}(\lambda)\mathcal{T}_{T-m}(\lambda)\,d\rho(\lambda),
\quad l,m=0,\ldots,T-1,\notag
\end{eqnarray}
and the following dynamic representation valid if $T\leqslant N$:
\begin{equation}
\label{MikhaylovAS_SP_mes_d1} C^T=C^T_N=
\begin{pmatrix}
r_0+r_2+\ldots+r_{2T-2} & r_1+\ldots+r_{2T-3} & \ldots &
r_T+r_{T-2} &
r_{T-1}\\
r_1+r_3+\ldots+r_{2T-3} & r_0+\ldots+r_{2T-4} & \ldots & \ldots
&r_{T-2}\\
\cdot & \cdot & \cdot & \cdot & \cdot \\
r_{T-3}+r_{T-1}+r_{T+1} &\ldots & r_0+r_2+r_4 & r_1+r_3 & r_2\\
r_{T}+r_{T-2}&\ldots &r_1+r_3&r_0+r_2&r_1 \\
r_{T-1}& r_{T-2}& \ldots & r_1 &r_0
\end{pmatrix}
\end{equation}
\end{proposition}

According to \cite{AT} the spectral measure $d\rho(\lambda)$
corresponding to operator $A$ give rise to the Fourier transform
$F: l^2\mapsto L_2(\mathbb{R},d\rho)$, defined by the rule:
\begin{equation}
\label{JM_Four} (Fa)(\lambda)=\sum_{n=0}^\infty
a_k\phi_k(\lambda),\quad a=(a_0,a_1,\ldots)\in l^2,
\end{equation}
where $\phi$ is a solution to (\ref{MikhaylovAS_Phi_def}). The
inverse transform and Parseval identity reads:
\begin{eqnarray}
a_k=\int_{-\infty}^\infty
(Fa)(\lambda)\phi_k(\lambda)\,d\rho(\lambda),\notag\\
\sum_{k=0}^\infty a_kb_k=\int_{-\infty}^\infty
{(Fa)(\lambda)}(Fb)(\lambda)\,d\rho(\lambda).\label{JM_parseval}
\end{eqnarray}
Note \cite{MikhaylovAS_MM2} that for $\lambda\in \mathbb{C}$ we
have the following representation for the Fourier transform of the
solution to (\ref{MikhaylovAS_Jacobi_dyn_int}) at $t=T$:
\begin{equation}
\label{JM_F_trans1}
\left(Fv^f_{\cdot,T}\right)(\lambda)=\sum_{k=1}^T
\mathcal{T}_k(\lambda)f_{T-k},\quad \lambda\in \mathbb{C}.
\end{equation}
Now we assume that $T=N$ and introduce the linear manifold of
Fourier images of states of dynamical system
(\ref{MikhaylovAS_Jacobi_dyn_int}) at time $t=N$, i.e. the Fourier
image of the reachable set:
\begin{equation*}
B_J^N:=F\mathcal{U}^N=\left\{\left(Fu^f_{\cdot,N}\right)(\lambda)\,|\,
f\in
\mathcal{F}^N\right\}=\left\{\left(Fv^f_{\cdot,N}\right)(\lambda)\,|\,
f\in \mathcal{F}^N\right\},
\end{equation*}

We equip $B_J^N$ with the scalar product defined by the rule:
\begin{gather}
[F,G]_{B^N_J}=\left(C^Nf,g\right)_{\mathcal{F}^N}, \quad F,G\in
B^N_J,\label{JM_scalprod}\\
F(\lambda)=\sum_{k=1}^N
\mathcal{T}_k(\lambda)f_{N-k},\,G(\lambda)=\sum_{k=1}^N
\mathcal{T}_k(\lambda)g_{N-k},\,f,g\in \mathcal{F}^N.\notag
\end{gather}
Evaluating (\ref{JM_scalprod}) making use of (\ref{JM_parseval})
yields:
\begin{eqnarray}
[F,G]_{B^N_J}=\left(v^f_{\cdot,N},v^g_{\cdot,T}\right)_{\mathcal{H}^N}=\left(u^f_{\cdot,N},u^g_{\cdot,T}\right)_{\mathcal{H}^N}
=\int_{-\infty}^\infty
{(Fu^f_{\cdot,N})(\lambda)}{(Fu^g_{\cdot,N})(\lambda)}\,d\rho(\lambda)\notag\\
=\int_{-\infty}^\infty
{F(\lambda)}\overline{G(\lambda)}\,d\rho(\lambda)=\int_{-\infty}^\infty
{F(\lambda)}\overline{G(\lambda)}\,d\rho_N(\lambda).\label{DBr_Scal}
\end{eqnarray}
On comparing (\ref{Moment_scal}) and (\ref{DBr_Scal}) we see that:
\begin{equation}
\label{Scal_DB_X} [F,G]_{B^N_J}=\langle
F,G\rangle=\sum_{n,m=0}^{N-1}s_{n+m}\alpha_n\overline{\beta_m}=\int_{-\infty}^\infty
{F(\lambda)}\overline{G(\lambda)}\,d\rho(\lambda).
\end{equation}
In \cite{MikhaylovAS_MM3} the authors proved the following
\begin{theorem}
\label{MikhaylovAS_Th_char} The vector
$(r_0,r_1,r_2,\ldots,r_{2N-2})$ is a response vector for the
dynamical system (\ref{MikhaylovAS_Jacobi_dyn_int}) if and only if
the matrix $C^T$ (with $T=N$) defined by
(\ref{MikhaylovAS_SP_mes_d}), (\ref{MikhaylovAS_SP_mes_d1}) is
positive definite.
\end{theorem}
This theorem shows that (\ref{DBr_Scal}) is a scalar product in
$B^N_J$. But we can say even more \cite{MikhaylovAS_MM2}:
\begin{theorem}
By dynamical system with discrete time  (\ref{Jacobi_dyn}) one can
construct the de Branges space
\begin{equation*}
B_J^N:=\left\{\left(Fu^f_{\cdot,N}\right)(\lambda)\,|\, f\in
\mathcal{F}^N\right\}=\left\{\sum_{k=1}^N
\mathcal{T}_k(\lambda)f_{N-k}\,|\, f\in \mathcal{F}^N\right\}.
\end{equation*}
As a set of functions it coincides with the space of Fourier
images of states of dynamical system (\ref{Jacobi_dyn}) at time
$N$ (the Fourier image of a reachable set) and is the set of
polynomials with real coefficients of the order less or equal to
$N-1$. The norm in $B_J^N$ is defined via the connecting operator:
\begin{equation*}
[F,G]_{B^N_J}:=\left(C^Nf,g\right)_{\mathcal{F}^N},\quad F,G\in
B^N_J,
\end{equation*}
where
\begin{equation*}
F(\lambda)=\sum_{k=1}^N
\mathcal{T}_k(\lambda)f_{N-k},\,G(\lambda)=\sum_{k=1}^N
\mathcal{T}_k(\lambda)g_{N-k},\,f,g\in \mathcal{F}^N.
\end{equation*}
The reproducing kernel has a form
\begin{equation*}
J_z(\lambda)=\sum_{k=1}^N \mathcal{T}_k(\lambda)j^z_{N-k},
\end{equation*}
where $j_z$ is a solution to Krein-type equation
\begin{equation*}
C^Nj^z=\overline{\begin{pmatrix}\mathcal{T}_{N}(z)\\
\mathcal{T}_{N-1}(z)\\\cdot\\ \mathcal{T}_1(z)\end{pmatrix}}.
\end{equation*}
\end{theorem}
Note \cite{MikhaylovAS_MM2,MM4} that control $j^z$ drives the
system (\ref{Jacobi_dyn}) to special state $\phi$, that is:
\begin{equation}
\label{Contr_eq}
\left(W^Nj^z\right)_i=\left(W^N_Nj^z\right)_i=\phi_i(z),\quad
i=1,\ldots,N.
\end{equation}
Thus the reproducing kernel in $C_N[X]$ (or in $B^N_J$) is given
by
\begin{equation}
\label{Rep_ker_CX}
K_N(z,\lambda)=\left(\left(C^N\right)^{-1}\overline{\begin{pmatrix}\mathcal{T}_{N}(z)\\
\mathcal{T}_{N-1}(z)\\\cdot\\
\mathcal{T}_1(z)\end{pmatrix}},{\begin{pmatrix}\mathcal{T}_{N}(\lambda)\\
\mathcal{T}_{N-1}(\lambda)\\\cdot\\
\mathcal{T}_1(\lambda)\end{pmatrix}}\right).
\end{equation}
\begin{remark}
The space of complex polynomials $C_N[X]$ with scalar product
defined by matrix $S$ is a de Branges space $B^N_J$ where scalar
product and reproducing kernel are given by (\ref{Scal_DB_X}) and
(\ref{Rep_ker_CX}).
\end{remark}
\begin{definition}
The $n-$th Christoffel function is defined by the rule
\begin{equation*}
\varkappa_n(\lambda)=\left(\sum_{k=1}^N\phi_k^2(\lambda)\right)^{-1}.
\end{equation*}
\end{definition}
From (\ref{Contr_eq}) it immediately follows that
\begin{equation*}
\varkappa_n(\lambda)=K_N(\lambda,\lambda)=\left(\left(C^N\right)^{-1}\overline{\begin{pmatrix}\mathcal{T}_{N}(\lambda)\\
\mathcal{T}_{N-1}(\lambda)\\\cdot\\
\mathcal{T}_1(\lambda)\end{pmatrix}},{\begin{pmatrix}\mathcal{T}_{N}(\lambda)\\
\mathcal{T}_{N-1}(\lambda)\\\cdot\\
\mathcal{T}_1(\lambda)\end{pmatrix}}\right).
\end{equation*}
Different formulas for reproducing kernel and Christoffel
functions are derived in \cite{MikhaylovAS_S,Schm}.

\section{Truncated moment problem. Recovering Dirichlet spectral data. }

We observe the following: in the moment problem we are given the
sequence of moments (\ref{MikhaylovAS_Moment_eq}), and in the
inverse dynamic problem for systems
(\ref{MikhaylovAS_Jacobi_dyn_int}), (\ref{Jacobi_dyn}) we are
given a response vector \cite{MikhaylovAS_MM,MikhaylovAS_MM3},
whose spectral representation has a form
(\ref{MikhaylovAS_Jac_resp_spectr}). Thus the knowledge of moments
$\{s_0,s_1,\ldots\}$ implies a possibility to calculate the
response vector $\{r_0,r_1,\ldots\}$ by
(\ref{MikhaylovAS_Jac_resp_spectr}). Note that Chebyshev
polynomials of the second kind $\{
\mathcal{T}_1(\lambda),\mathcal{T}_2(\lambda),\ldots
\mathcal{T}_n(\lambda)\}$ (see (\ref{Cheb})) are related to $\{
1,\lambda,\lambda^{n-1}\}$ by the following formula
\begin{equation}
\label{Response_moments_rel}
\begin{pmatrix}
\mathcal{T}_1(\lambda)\\
\mathcal{T}_2(\lambda)\\
\ldots \\
\mathcal{T}_n(\lambda)
\end{pmatrix}=\Lambda_n \begin{pmatrix}
1\\
\lambda\\
\ldots \\
\lambda^{n-1}
\end{pmatrix}=\begin{pmatrix}
1& 0& \ldots & 0\\
a_{21}& 1 & \ldots & 0\\
\ldots \\
a_{n1} & a_{n2}& \ldots & 1
\end{pmatrix}\begin{pmatrix}
1\\
\lambda\\
\ldots \\
\lambda^{n-1}
\end{pmatrix}.
\end{equation}
\begin{proposition}
Entries of the matrix $\Lambda_n\in R^{n\times n}$ are given by
\begin{equation}
\label{LambdaMatr}
\Lambda_n=a_{ij}=\begin{cases} 0,\quad \text{if $i>j$},\\
0,\quad \text{if $i+j$ is odd,}\\
C_{\frac{i+j}{2}}^j(-1)^{\frac{i+j}{2}+j}.
\end{cases}
\end{equation}
entries of the response vector are related to moments by the rule:
\begin{equation}
\label{Resp_and_moments}
\begin{pmatrix}
r_0\\
r_1\\
\ldots \\
r_{n-1}
\end{pmatrix}=\Lambda_n\begin{pmatrix}
s_0\\
s_1\\
\ldots \\
s_{n-1}
\end{pmatrix}.
\end{equation}
\end{proposition}

\begin{definition}
By a solution of a truncated moment problem of order $N$ we call a
Borel measure $d\tilde\rho_N(\lambda)$ on $\mathbb{R}$ such that
equalities (\ref{MikhaylovAS_Moment_eq}) with this measure hold
for $k=0,1,\ldots,2N.$
\end{definition}

\begin{remark}
Results from \cite{MikhaylovAS_MM,MikhaylovAS_MM3} imply that from
the finite set of moments $\{s_0,s_1,\ldots,s_{2N-2}\}$, or what
is equivalent from $\{r_0,r_1,\ldots,r_{2N-2}\}$, it is possible
to recover Jacobi matrix $A^N\in \mathbb{R}^{N\times N}$, whose
elements can be thought of as a coefficients in dynamical system
(\ref{MikhaylovAS_Jacobi_dyn_int}) with Dirichlet boundary
condition at $n=N+1$, or $N\times N$ block in semi-infinite Jacobi
matrix in (\ref{MikhaylovAS_Jacobi_dyn_int}).
\end{remark}

This theorem and formulas for the entries of Jacobi matrix
obtained in \cite{MikhaylovAS_MM3} implies the following procedure
of solving the truncated moment problem:
\begin{itemize}
\item{1)} Calculate $(r_0,r_1,r_2,\ldots,r_{2N-2})$ from
$\{s_0,s_1,\ldots,s_{2N-2}\}$ by using (\ref{Resp_and_moments}).

\item{2)} Recover $N\times N$ Jacobi matrix $A^N$ using formulas
for $a_k,$ $b_k$ from \cite{MikhaylovAS_MM3}

\item{3)} Recover spectral measure for finite Jacobi matrix $A^N$
prescribing arbitrary selfadjoint condition at $n=N+1$. Or one can
do

\item{3')} Extend Jacobi matrix $A^N$ to \emph{finite} Jacobi
matrix $A^M$, $M>N$, prescribe arbitrary selfadjoint condition at
$n=M+1$ and recover the spectral measure of $A^M$.

\item{3'')} Extend Jacobi matrix $A^N$ to \emph{infinite} Jacobi
matrix $A$, and recover the spectral measure of $A$.
\end{itemize}
Every measure obtained in $3),$ $3')$, $3'')$ gives a solution to
the truncated moment problem. Below we propose a different
approach: we recover the spectral measure corresponding to Jacobi
matrix directly from moments (from the operator $C^N$), without
recovering the Jacoi matrix itself.

\begin{agreement}
\label{MikhaylovAS_argr} We assume that controls $f\in
\mathcal{F}^N$, $f=\left(f_0,\ldots,f_{N-1}\right)$ are extended:
$f=\left(f_{-1},f_0,\ldots,f_{N-1},f_N\right)$, where $f_1=f_N=0$.
\end{agreement}
We introduce the special space of controls
$\mathcal{F}^N_0=\left\{f\in \mathcal{F}^T\,|\, f_0=0\right\}$ and
the operator $D: \mathcal{F}^T\mapsto \mathcal{F}^T$ acting by
\begin{equation*}
\left(Df\right)_t=f_{t+1}+f_{t-1}.
\end{equation*}
The following statements can be easily proved using arguments from
\cite{MikhaylovAS_MM3}:
\begin{proposition}
\label{MikhaylovAS_PropContr} The operator $W^N$ maps
$\mathcal{F}^N$ isomorphically onto $\mathcal{H}^{N}$ and
$\mathcal{F}^N_0$ maps isomorphically onto $\mathcal{H}^{N-1}$.
\end{proposition}
\begin{proposition}
On the set $\mathcal{F}^T_0$ the following relation holds:
\begin{equation}
\label{MikhaylovAS_FT_proper} W^NDf=DW^Nf,\quad f\in
\mathcal{F}^N_0.
\end{equation}
\end{proposition}
Taking $f,g\in \mathcal{F}_0^T$ we can evaluate the quadratic
form, bearing in mind (\ref{MikhaylovAS_FT_proper}):
\begin{equation}
\label{MikhaylovAS_L_quadr}
\left(C^NDf,g\right)_{\mathcal{F}^N}=\left(W^NDf,W^Ng\right)_{\mathcal{H}^{N}}=
\left(DW^Nf,W^Ng\right)_{\mathcal{H}^{N-1}}=\left(A^{N-1}v^f,v^g\right)_{\mathcal{H}^{N-1}}.
\end{equation}
The last equality in (\ref{MikhaylovAS_L_quadr}) means that only
$A^{N-1}$ block from the whole matrix $A^N$ is in use. Then it is
possible to perform the spectral analysis of $A^{N-1}$ using the
classical variational approach, the controllability of the system
(\ref{MikhaylovAS_Jacobi_dyn_int}) (see Proposition
\ref{MikhaylovAS_PropContr}) and the representation
(\ref{MikhaylovAS_L_quadr}), see also \cite{MikhaylovAS_B2001}.
The spectral data of Jacobi matrix $A^{N-1}$ with Dirichlet
boundary condition at $n=N$ can be recovered by the following
procedure:
\begin{itemize}
\item[1)] The first eigenvalue is given by
\begin{equation}
\label{MikhaylovAS_M1}\lambda_1^{N-1}=\min_{f\in
\mathcal{F}^N_0,\,(C^N f,f)_{\mathcal{F}^N}=1}\left(C^N
Df,f\right)_{\mathcal{F}^N}.
\end{equation}

\item[2)] Let $f^1$, be the minimizer of (\ref{MikhaylovAS_M1}),
then
\begin{equation*}
\rho_1=\left(C^N f^1,f^1\right)_{\mathcal{F}^N}.
\end{equation*}

\item[3)] The second eigenvalue is given by
\begin{equation}
\label{MikhaylovAS_M2}\lambda_2^{N-1}=\min_{\substack {f\in
\mathcal{F}^N_0,(C^N f,f)_{\mathcal{F}^N}=1\\
(C^N f,f_l)_{\mathcal{F}^N}=0}}\left(C^N
Df,f\right)_{\mathcal{F}^N}.
\end{equation}

\item[4)] Let $f^2$, be the minimizer of (\ref{MikhaylovAS_M2}),
then
\begin{equation*}
\rho_2=\left(C^T f^2,f^2\right)_{\mathcal{F}^T}.
\end{equation*}
\end{itemize}
Continuing this procedure, one recovers the set
$\{\lambda_k^{N-1},\rho_k\}_{k=1}^{N-1}$ and construct the measure
$d\rho^{N-1}(\lambda)$ by formula (\ref{MikhaylovAS_measure}).

\begin{remark}
The measure, constructed by the above procedure solves the
truncated moment problem for the set of moments
$\{s_0,s_1,\ldots,s_{2N-4}\}$.
\end{remark}

By $f^k$, $k=1,\ldots,N$ we denote the control that drive system
(\ref{MikhaylovAS_Jacobi_dyn_int}) to prescribed state (see
(\ref{MikhaylovAS_Phi_def})):
\begin{equation*}
W^Tf_k=\phi^k,\quad k=1,\ldots,N.
\end{equation*}
Due to Proposition \ref{MikhaylovAS_PropContr}, such a control
exists and is unique for every $k$. The remarkable fact that these
controls as well as the spectrum of $A^N$ can be found from
Euler-Lagrange equations for the problem of the minimization of a
functional $\left(C^NDf,f\right)_{\mathcal{F}^N}$ in
$\mathcal{F}^N_0$ with the constrain
$\left(C^Tf,f\right)_{\mathcal{F}^N}=1$. Similar method of
deriving equations which can be used for recovering of spectral
data was used in \cite{MikhaylovAS_AMM}. We introduce the operator
\begin{equation}
\label{BN_def}
B^N=\begin{pmatrix} c_{N,N+1}+c_{N,N-1} & c_{N,N}+c_{N,N-2} & \ldots & c_{N,3}+c_{N,1} & c_{N,2}\\
c_{N-1,N+1}+c_{N-1,N-1} & \ldots & \ldots & c_{N-1,3}+c_{N-1,1} & c_{N-1,2}\\
\cdot & \cdot & \ldots & \cdot & \cdot \\
c_{1,N+1}+c_{1,N-1}& c_{1,N}+c_{1,N-2} & \ldots & \ldots & c_{1,2}
\end{pmatrix}.
\end{equation}
The following result was obtained in  \cite{MM4}:
\begin{theorem}
\label{MikhaylovAS_teor} The spectrum of $A^{N}$ and
(non-normalized) controls $f_k$, $k=1,\ldots,N$ are the spectrum
and eigenvectors of the following generalized spectral problem:
\begin{equation}
\label{MikhaylovAS_m_eqn} B^Nf_k=\lambda_kC^Nf_k,\quad
k=1,\ldots,N.
\end{equation}
\end{theorem}
Introduce the following Hankel matrices
\begin{equation*}
S^N_{m}:=\begin{pmatrix} s_{2N-2+m} & s_{2N-3+m} & \ldots & s_{N-1+m}\\
s_{2N-3+m} & \ldots & \ldots & \ldots\\
\cdot & \cdot & \ldots & s_{1+m} \\
s_{N-1+m} & \ldots & s_{1+m} & s_{m}
\end{pmatrix},\quad m=0,1,\dots,
\end{equation*}
the matrix $J_N\in \mathbb{R}^{N\times N}$:
\begin{equation*}
J_N=\begin{pmatrix} 0 & \ldots & 0 & 1\\
0 & \ldots & 1 & 0\\
\cdot & \cdot & \ldots & \cdot \\
0& 1 & \ldots & 0\\
 1& \ldots & 0 & 0
\end{pmatrix},\quad J_NJ_N=I_N=\begin{pmatrix} 1 & 0 & \ldots & 0\\
0 & 1 & \ldots & 0\\
\cdot & \cdot & \ldots & \cdot \\
0& \ldots & 1 & 0\\
0& \ldots & 0 & 1
\end{pmatrix},
\end{equation*}
and define
\begin{equation*}
\widetilde\Lambda_N:=J_N\Lambda_N J_N.
\end{equation*}
The remarkable fact is that the matrices $B^N,\,C^N$ can be
reduced to Hankel matrices by the same linear transformation:
\begin{theorem}
\label{Propos_EL} The following relations hold:
\begin{eqnarray*}
C^N=\widetilde\Lambda_N S_0^N
\left(\widetilde\Lambda_N\right)^*,\\
B^N =\widetilde\Lambda_N S_1^N \left(\widetilde\Lambda_N\right)^*.
\end{eqnarray*}
Then the generalized spectral problem (\ref{MikhaylovAS_m_eqn})
upon introducing the notation
$g_k=\left(\widetilde\Lambda_N\right)^*f_k$ is equivalent to the
following generalized spectral problem:
\begin{equation}
\label{m_eqn2} S^N_1g_k=\lambda_k S^N_0g_k.
\end{equation}
\end{theorem}
Having found spectrum and non-nomalized controls from
(\ref{MikhaylovAS_m_eqn}) one can recover the measure of $A^N$
with Dirichlet boundary condition at $n=N+1$ by the following
procedure:
\begin{itemize}

\item{1)} Normalize controls by choosing
$\left(C^Nf_k,f_k\right)_{\mathcal{F}^N}=1$,

\item{2)} Observe that $W^Nf^k=\alpha_k\phi^k$ for some
$\alpha_k\in \mathbb{R}$, where the constant is defined by
$\alpha_k=\left(Rf_k\right)_N$.

\item{3)} The norming coefficients are given by
$\rho_k={\alpha_k^2}$, $k=1,\ldots,N$.

\item{4)} Recover the measure by (\ref{MikhaylovAS_measure}).

\end{itemize}

\section{Existence and uniqueness for Hamburger, Stieltjes and Hausdorff moment
problems.}

We remind the reader that the moment problem is called
\emph{determinate} if it has only one solution, otherwise it is
called \emph{indeterminate}. It is well-known fact
\cite{MikhaylovAS_S,Schm} that the uniqueness of the solution to a
moment problem is related to the index of the operator $A$. Here
we provide well-known results on discrete version of Weyl limit
point-circle theory, answering the question on the index of $A$,
which will be subsequently used. By $\xi(\lambda)$ we denote the
solution to the difference equation in (\ref{MikhaylovAS_Phi_def})
with Cauchy data $\xi_0=-1$, $\xi_1=0$.

\begin{proposition}
\label{Hamburher_crt} The Jacobi operator $A$ is limit circle at
infinity (has index equal to one) if and only if one of the
following occurs:
\begin{itemize}


\item[1)] $\phi(x),\xi(x)\in l^2$ for some $x\in \mathbb{R}$,


\item[2)] $\phi(x),\varphi'(x)\in l^2$ for some $x\in \mathbb{R}$,

\item[3)] $\xi(x),\xi'(x)\in l^2$ for some $x\in \mathbb{R}$.
\end{itemize}
\end{proposition}
In \cite{MM4} the authors proved the following
\begin{theorem}
\label{Th_char_new} The  set of numbers $(s_0,s_1,s_2,\ldots)$ are
moments of a spectral measure corresponding to the Jacobi operator
$A$ if and only if
\begin{equation}
\label{HamCnd} \text{the matrix $S^N_0$ is positive definite for
all $N\in \mathbb{N}$.}
\end{equation}
The Hamburger moment problem is indeterminate if and only if
\begin{equation}\label{Ham_cond}
\lim_{T\to\infty}\left(\left(C^T\right)^{-1}\Gamma_T,\Gamma_T\right)_{\mathcal{F}^T}<+\infty,\quad
\lim_{T\to\infty}\left(\left(C^T\right)^{-1}\Delta_T,\Delta_T\right)_{\mathcal{F}^T}<+\infty,
\end{equation}
where
\begin{equation} \label{GaOm_def}
\Gamma_T:=\begin{pmatrix} \mathcal{T}_T(0)\\
\mathcal{T}_{T-1}(0)\\\ldots\\
\mathcal{T}_1(0)\end{pmatrix},\quad \Omega_T=\begin{pmatrix} \mathcal{T}_T'(0)\\
\mathcal{T}_{T-1}'(0)\\\ldots\\
\mathcal{T}_1'(0)\end{pmatrix}.
\end{equation}
\end{theorem}
Here we rewrite conditions (\ref{Ham_cond}) in more standard form:
\begin{proposition}
Conditions in (\ref{Ham_cond}) are equivalent to
    \begin{equation}
    \label{Ham_cond_S} \lim_{T\to\infty} \frac{\det S_2^{T-1}}{\det S^T_0}    <+\infty,\quad
    \lim_{T\to\infty} \frac{\det S_0^{T-1,2}}{\det S^T_0} <+\infty,
    \end{equation}
    \begin{equation*}
    \text { where}\quad S_0^{T-1,2}=\begin{pmatrix}
        s_0& s_2& \ldots & s_{T-2}\\
        s_{2}& s_{4} & \ldots & s_{T-1}\\
        \ldots \\
        s_{T-2} & s_{T-1}& \ldots & s_{2T-3}
    \end{pmatrix}.
    \end{equation*}
\end{proposition}
 Indeed, bearing in mind (\ref{Response_moments_rel}) and relations
 $$
 C^T=\tilde\Lambda_T S_0^T \tilde\Lambda_T^*,\quad \tilde\Lambda_T=J_T \Lambda_T
 J_T,
 $$
 we pass to
\begin{eqnarray}
 \left(\left(C^T\right)^{-1}\Gamma_T,\Gamma_T\right)_{\mathcal{F}^T}= \left(\left(\tilde
 S^T_0\right)^{-1}e_1,e_1\right)_{\mathcal{F}^T},\label{SS_1}\\
  \left(\left(C^T\right)^{-1}\Delta_T,\Delta_T\right)_{\mathcal{F}^T}= \left(\left(\tilde S^T_0\right)^{-1}e_2,e_2\right)_{\mathcal{F}^T}.\label{SS_2}
\end{eqnarray}
The right hand side of above equalities can be
computed using the following formula for a bilinear form of
inverse matrix. Namely, for the matrix $D=(d_{ij})_{i,j=1}^n$ and
vectors $h=(h_1,\dots,h_n)$, $c=(c_1,\dots,c_n)$ we have that
 \begin{equation}\label{inv_bf}
    (D^{-1}b,c)=\frac{\det{\begin{pmatrix}
            0& h_1& \ldots & h_n\\
            c_{1}& d_{11} & \ldots & d_{1n}\\
            \ldots \\
            c_{n} & d_{n1}& \ldots & d_{nn}
            \end{pmatrix}}}{\det{\begin{pmatrix}

                 d_{11} & \ldots & d_{1n}\\
                \ldots \\
                 d_{n1}& \ldots & d_{nn}
                \end{pmatrix}}}.
 \end{equation}
Then applying (\ref{inv_bf}) to (\ref{SS_1}), (\ref{SS_2}), we get
(\ref{Ham_cond_S}).

The following result concerning the Stieltjes problem was
formulated in \cite{MM4}, where the authors obtained expressions
for the mass $M_\infty$ and the length $L_\infty$ of a string in
terms of operators of the Boundary control method, associated with
the dynamical system (\ref{Jacobi_dyn}):
\begin{theorem}
The  set of numbers $(s_0,s_1,s_2,\ldots)$ are moments of a
spectral measure, supported on $(0,+\infty)$, corresponding to
Jacobi operator $A$
\begin{equation*}
\text{matrices $S^N_0$ and $S^N_1$ are positive definite for all
$N\in \mathbb{N}$.}
\end{equation*}
The Stieltjes moment problem is indeterminate if and only if the
following relations hold:
\begin{equation}
\label{Stiel_cond}
M_\infty=\lim_{T\to\infty}\left(\left(C^T\right)^{-1}\Gamma_T,\Gamma_T\right)_{\mathcal{F}^T}<+\infty,\quad
L_\infty=\lim_{K\to\infty}\frac{\left(\left(C^{K}\right)^{-1}\left(R^{K}\right)^*\Gamma_{K},e_1\right)}
{\left(\left(C^{K}\right)^{-1}\Gamma_{K},e_1\right)}<+\infty.
\end{equation}
\end{theorem}
Notice that the the necessity of (\ref{Stiel_cond}) is a subtle
result, see  (see \cite{MikhaylovAS_S,Schm}) for details. Here we
reformulate (\ref{Stiel_cond}):
\begin{proposition}
Conditions in (\ref{Stiel_cond}) are equivalent to
\begin{equation}
\label{Stiel_cond_S} \lim_{T\to\infty} \frac{\det S_2^{T-1}}{\det S^T_0}    <+\infty,\quad
\lim_{T\to\infty} \frac{\det S_{0,0}^{T}}{\det S^{T-1}_1} <+\infty,
\end{equation}
\begin{equation*}
 \text{ where }
S_{0,0}^{T}=\begin{pmatrix}
0& s_0& \ldots & s_{T-2}\\
s_{0}& s_{1} & \ldots & s_{T-1}\\
\ldots \\
s_{T-1} & s_{T}& \ldots & s_{2T-2}
\end{pmatrix}.
\end{equation*}
\end{proposition}
 Indeed, we notice that the operator $R^T$ has a form of a Toeplitz matrix
 $$
   R^T=\begin{pmatrix}
   0& 0& \ldots & 0\\
   r_0& 0 & \ldots & 0\\
   \ldots \\
   r_{T-2} & r_{T-3}& \ldots & 0
   \end{pmatrix}.
 $$
Bearing in mind (\ref{Response_moments_rel}) and relations
$$
C^T=\tilde\Lambda_T S_0^T \tilde\Lambda_T^*,\quad
\tilde\Lambda_T=J_T \Lambda_T J_T,
$$
we pass to
$$
\left(\left(C^T\right)^{-1}\Gamma_T,e_1\right)_{\mathcal{F}^T}= \left(\left(\tilde S^T_0\right)^{-1}e_1,e_T\right)_{\mathcal{F}^T},
$$
$$
\left(\left(C^T\right)^{-1}(R^T)^* \Gamma_T,e_1\right)_{\mathcal{F}^T}= \left(\left(\tilde S^T_0\right)^{-1}g,e_T\right)_{\mathcal{F}^T}.
$$
where
$g=\Lambda_T^{-1}J_T(R^T)^*\Gamma_T=(0,s_0,s_1,\dots,s_{T-2})$.
Using (\ref{inv_bf}) we obtain (\ref{Stiel_cond_S}).

\begin{remark}
Condition $\xi(0)\in l^2$ is equivalent to
    \begin{equation*}
    \left(\left(C^T\right)^{-1} \left(R^T\right)^*\Gamma_T,\left(R^T\right)^*\Gamma_T\right)_{\mathcal{F}^T}=
    \left(\left(\tilde S^T_0\right)^{-1} g,g\right)_{\mathcal{F}^T} <\infty,
    \end{equation*}
    where $g=(0,s_0,s_1,\dots,s_{T-2})$.
\end{remark}

The following stetement was proved in \cite{MikhaylovAS_S}, but we
give a new proof in order to show that it is a direct consequence
of the spectral problem (\ref{MikhaylovAS_m_eqn}).

\begin{proposition}
Let $\{s_k\}_{k=0}^\infty$ be a set of Stieltjes moments, we set
    $\{h_m\}_{m=0}^\infty=\{s_0,0,s_1,0,s_2,0,s_3,0,\dots\}$. Then
    $\{h_m\}_{m=0}^\infty$ corresponds to a determinate Hamburger moment problem if and only if
    $\{s_k\}_{k=0}^\infty$ corresponds to a determinate Stieltjes moment problem.
\end{proposition}
We use the spectral problem (\ref{m_eqn2}) with two matrices for
Stiltjes moment problem:
\begin{equation}
\label{St_matr} \tilde S_0^T=\begin{pmatrix}
s_0& s_1& \ldots & s_{T-1}\\
s_1& s_2 & \ldots & s_{T}\\
\ldots \\
s_{T-1} & s_{T}& \ldots & s_{2T-2}
\end{pmatrix},\quad
 \tilde S_1^T=\begin{pmatrix}
s_1& s_2& \ldots & s_{T}\\
s_2& s_3 & \ldots & s_{T+1}\\
\ldots \\
s_{T} & s_{T+1}& \ldots & s_{2T-1}
\end{pmatrix}
\end{equation}
and two matrices for Hamburger moment problem:
\begin{equation}
\label{Ham_matr} \tilde H_0^{2T}=\begin{pmatrix}
s_0&   0& s_1& \ldots & 0\\
  0& s_1&   0& \ldots & s_{T}\\
s_1&   0&   s_2& \ldots & 0\\
\ldots \\
 0& s_{T}& 0&  \ldots & s_{2T-1}
\end{pmatrix},\quad
\tilde H_1^{2T}=\begin{pmatrix}
0&   s_1& 0& \ldots & s_T\\
s_1& 0&   s_2& \ldots & 0\\
0&   s_2&   0& \ldots & s_{T+1}\\
\ldots \\
s_T& 0& s_{T+1}&  \ldots & 0
\end{pmatrix}
\end{equation}
Simple observation $\det{\tilde H_0^{2T}}  = \det{\tilde S_0^{T}}
\det{\tilde S_1^{T}}$ allows us to conclude that the set
$\{h_m\}_{m=0}^\infty=\{s_0,0,s_1,0,s_2,0,s_3,0,\dots\}$ indeed
corresponds to Hamburger moment problem. Even more, if we look for
eigenvalues of (\ref{m_eqn2}) with Hamburger matrices
(\ref{Ham_matr}), we see that
\begin{multline*}
0=\det{(\lambda \tilde H_0^{2T} -\tilde
H_1^{2T})}=\det{\begin{pmatrix}
    \lambda s_0&   -s_1& \lambda s_1& \ldots & -s_T\\
    -s_1& \lambda s_1&   -s_{2}& \ldots & \lambda s_{T}\\
    \lambda s_1&   -s_{2}&   \lambda s_2& \ldots & -s_{T+1}\\
    \ldots \\
    -s_{T}& \lambda s_{T}& -s_{T+1}&  \ldots & \lambda s_{2T-1}
    \end{pmatrix}}=\\
=\det{\begin{pmatrix}
    \lambda^2 s_0-s_1&   0& \lambda^2 s_1-s_2& \ldots & 0\\
    0&  s_1&   0& \ldots & \lambda^2 s_{T}-s_{T+1}\\
    \lambda^2 s_1-s_2&   0&   \lambda^2 s_2-s_3& \ldots & 0\\
    \ldots \\
    0& s_{T}& 0&  \ldots &  s_{2T-1}
    \end{pmatrix}} = \det{(\lambda^2 S_0^T - S_1^T)}\det{S_1^T}.
\end{multline*}
The later means that if $0<\mu$ is an eigenvalue for
(\ref{m_eqn2}) with Stieltjes matrices (\ref{St_matr}), then
$\pm\sqrt{\mu}$\,--\, are two eigenvalues for (\ref{m_eqn2}) with
Hamburger matrices (\ref{Ham_matr}) and vise-versa.

To prove that Hamburger and Stieltjes problems are determinate simultaneously we note that
$$
\frac{\det{H_2^{2T-1}}}{\det{H_0^{2T}}} =\frac{\det{S_1^{T}}\det{S_2^{T-1}}}{\det{S_0^{T}}\det{S_1^{T}}} =
\frac{\det{S_2^{T-1}}}{\det{S_0^{T}}}
$$
and
$$
\left(\left(\tilde H^{2T}_0\right)^{-1} g,g\right)_{\mathcal{F}^{T}} = \frac{1}{\det{H_0^{2T}}}\det\begin{pmatrix}
0&   0&  s_0& \ldots & s_{T-1}\\
0&  s_0&   0& \ldots &  0\\
s_0&   0&    s_1& \ldots & s_{T}\\
\ldots \\
s_{T-1}& 0& s_{T}&  \ldots &  s_{T-1}
\end{pmatrix}=\frac{\det{S_0^{T}}\det{S_{0,0}^{T+1}}}{\det{S_0^{T}}\det{S_1^{T}}}=
\frac{\det{S_{0,0}^{T+1}}}{\det{S_1^{T}}} .
$$
Using Theorem 7 for Stieltjes problem and Theroem 5 and Remark 5
for Hamburger we complete the proof.

As we saw, the reproducing kernel $K^N(z,\lambda)$ has a form
(\ref{Rep_ker_CX}). Using (\ref{inv_bf}) and carrying out similar
transformations, we obtain the following
\begin{remark}
The reproducing kernel admits the following representation:
\begin{equation*}
K^N(z,\lambda) = \frac{1}{\det{S_0^{T}}} \det\begin{pmatrix}
0&   1&    z& \ldots & z^{N-1}\\
1&  s_0&   s_1& \ldots &  s_{T-1}\\
\lambda&  s_1&   s_2& \ldots &  s_{T}\\
\ldots \\
\lambda^{N-1}& s_{T-1}& s_{T}&  \ldots &  s_{2T-2}
\end{pmatrix}
\end{equation*}
\end{remark}

In \cite{MM4} the authors used (\ref{m_eqn2}) to prove the
following
\begin{theorem}
The  set of numbers $(s_0,s_1,s_2,\ldots)$ are moments of a
spectral measure, supported on $(0,1)$, corresponding to operator
$A$ if and only if the condition
\begin{equation*}
S^N_0\geqslant S^N_1 > 0\quad \text{holds for all $N\in
\mathbb{N}$}
\end{equation*}
The Hausdorff moment problem is determinate.
\end{theorem}

\subsection*{Conclusion}

In the present paper we considered the dynamical system
(\ref{MikhaylovAS_Jacobi_dyn_int}) associated with Jacobi matrix
with the Dirichlet boundary condition at the "right end". In the
third section it is shown how to use the Boundary control method
to recover the measure $d\rho_N$ of the operator $A^N$
(\ref{A_eq1}), (\ref{A_eq2}) from the finite set of moments. After
taking the limit $d\rho^N(\lambda)\mapsto d\rho^*(\lambda)$, where
convergence is understood in the week sense, we have two options:
when $A$ is in limit point case at infinity, $d\rho^*(\lambda)$ is
a unique solution to a moment problem, but when $A$ is in limit
circle case, $d\rho^*(\lambda)$ gives certain distinguished
solution (week limit of measures corresponding to operators with
Dirichlet condition). At the same time it is well-known
\cite{MikhaylovAS_Ahiez,MikhaylovAS_S,Schm} that the answer in
Hamburger moment problem is a spectral measures of any
self-adjoint extension of Jacobi operator $A$.

The prospective problem is to use dynamic inverse data
$\{r_0,r_1,\ldots\}$ of the dynamical system (\ref{Jacobi_dyn})
associated with $A$, which is obtained from the set of moments
$\{s_0,s_1,\ldots\}$ by (\ref{Resp_and_moments}), for the
construction of the dynamic model of self-adjoint extensions of
operator $A$ in the spirit of \cite{BS1,BS2}.

\subsection*{Acknowledgments}

The research of Victor Mikhaylov was supported in part by RFBR
17-01-00529, RFBR 18-01-00269 and the Ministry of Education and
Science of Republic of Kazakhstan under grant AP05136197. Alexandr
Mikhaylov was supported by RFBR 17-01-00099 and RFBR 18-01-00269.

\end{document}